\documentclass[11pt]{article}

\usepackage[margin=1in]{geometry}
\usepackage{orcidlink}
\usepackage{amsmath,amssymb,amsthm,enumerate}
\usepackage{enumitem}
\usepackage{hyperref}
\usepackage{float}
\usepackage[round]{natbib}
\usepackage{bbm}
\usepackage{tikz}
\usetikzlibrary{arrows.meta, positioning}
\newtheorem*{portmanteautheorem}{Theorem (Portmanteau)}
\newtheorem*{mct}{Theorem (Monotone Convergence Theorem)}
\newtheorem*{bct}{Theorem (Bounded Convergence Theorem)}
\newtheorem*{skorohod}{Theorem (Skorokhod's Representation Theorem)}

\newtheorem*{Proof}{Proof}
\newcommand{\ind}{\mathbbm{1}}

\begin{document}
\begin{center}
{\Large \textbf{A Self-Contained Proof of the Portmanteau Theorem}}\\[0.5em]
{\normalsize Benjamin Smith\orcidlink{0009-0007-2206-0177}}\\
{\normalsize Department of Statistical Sciences, University of Toronto}\\
{\normalsize $\texttt{benyamin.smith@mail.utoronto.ca}$}
\end{center}

\vspace{1em}

\noindent\textbf{Abstract.}
The Portmanteau Theorem gives several equivalent characterizations of weak convergence of probability measures (equivalently, convergence in distribution of random variables). This note presents a unified proof of the theorem in a seven-statement formulation, emphasizing the structure of the implication cycle connecting expectations of bounded continuous functions with probability bounds on open and closed sets. The argument highlights standard approximation techniques and clarifies the role of intermediate function classes such as bounded Lipschitz functions.

\section{Preliminaries}
 To keep this proof self contained, we will state Skorokhod's Representation Theorem, the Monotone Convergence Theorem and the Bounded Convergence Theorem as they will be applied in the proof. For brevity, their proofs will not be shown. The reader is invited to read works such as \cite{Billingsley1999}, \cite{Pollard_2001}, \cite{Rosenthal2006}, \cite{durrett:2019} and other classic graduate probability texts for a more rigorous treatment.

\begin{skorohod}
Let $F_n$ and $F$ be distribution functions.  If $F_n \Rightarrow F$, then there exists random variables $Y_n$ and $Y$ with distribution functions $F_n$ and $F$ respectively such that $Y_n \overset{\text{a.s.}}{\longrightarrow}Y$
\end{skorohod}

\begin{mct}
If $X_1, X_2, \dots$ are non-negative random variables and  $\{X_n\}\nearrow X$. Then
$$
\lim_{n \rightarrow\infty}\mathbb{E}[X_n] = \mathbb{E}[X]
$$
\end{mct}

\begin{bct}
Let $Z_n$ be a sequence of bounded random variables such that $|Z_n| \le M$ for some constant $M$. If $Z_n \xrightarrow{a.s.} Z$, then 
$$ \lim_{n \rightarrow \infty} \mathbb{E}[Z_n] = \mathbb{E}[Z] $$
\end{bct}
\bigskip

\section{The Portmanteau Theorem and Proof}
\subsection{Portmanteau Theorem}
To show all properties explicitly, we state the Portmanteau Theorem similarly to how it is presented by \citet{vanderVaart1998} with seven statements. For a five statement version that subsumes some of these properties, see \citep{Billingsley1999} or \citep{durrett:2019}.
\begin{portmanteautheorem}
Let $X_n$ and $X$ be random vectors. The following are equivalent
\begin{description}
    \item[(i)] $X_n \overset{d}{\Rightarrow} X$ ($X_n$ converges in distribution to $X$).
    \item[(ii)] $\forall$ bounded continuous functions $g$, $\underset{n \rightarrow \infty}\lim \mathbb{E}[g(X_n)] = \mathbb{E}[g(X)]$
    \item[(iii)] $\forall$ bounded continuous and Lipschitz\footnote{A Lipschitz function is a function whose rate of change is limited, meaning a small change in input leads to a proportionally small change in output, bounded by a constant \(K\) (the Lipschitz constant) such that $$|f(x) - f(y)| \le K|x - y|$$ 
    
    for all $x, y$ in the domain.} function $g$, $\underset{n \rightarrow \infty}\lim \mathbb{E}[g(X_n)] = \mathbb{E}[g(X)]$
    \item[(iv)] $\forall$ non-negative continuous functions $g$, $\lim\inf_{n \longrightarrow \infty}\mathbb{E}[g(X_n)] \ge \mathbb{E}[g(X)]$
    \item[(v)] $\forall$ open sets $G$ $\lim\inf_{n \rightarrow \infty}P(X_n \in G) \ge P(X \in G)$ 
    \item[(vi)] $\forall$ closed sets $K$ $\limsup_{n \rightarrow \infty}P(X_n \in K) \le P(X \in K)$
    \item[(vii)] $\forall$ Borel sets $A$ with $P(X\in \partial A) = 0$, $\lim_{n \rightarrow \infty}P(X_n \in A) = P(X\in A)$\\ ($\partial A$ denotes the boundary of A)
\end{description} 
\end{portmanteautheorem}
\noindent

\subsection{Proof}
\begin{figure}[h!]
    \centering
    \begin{tikzpicture}[
        node distance=1.2cm and 1cm,
        every node/.style={inner sep=2pt,font=\bfseries},
        thick,
        >={Stealth[length=2mm]} % Makes the arrowheads sharper
    ]
        % Define the main horizontal nodes
        \node (i)   {(i)};
        \node (ii)  [right=of i]   {(ii)};
        \node (iii) [right=of ii]  {(iii)};
        \node (v)   [right=of iii] {(v)};
        \node (vi)  [right=of v]   {(vi)};
        
        % Define the outlier nodes
        \node (iv)  [above right=0.6cm and 0.2cm of ii] {(iv)};
        \node (vii) [below right=0.8cm and 0.2cm of v]  {(vii)};

        % Draw the simple implications
        \draw[->] (i)   -- (ii);
        \draw[->] (ii)  -- (iii);
        \draw[->] (iii) -- (v);
        \draw[->] (v)   -- (vi);
        
        % Draw the double-headed arrow for (ii) <=> (iv)
        \draw[<->] (ii) -- (iv);
        
        % Draw the descending arrows to (vii)
        \draw[->] (v)  -- (vii);
        \draw[->] (vi) -- (vii);
        
        % Draw the long "Cycle" return arrow
        % Using "shorten" to ensure it doesn't overlap the text
        \draw[->, shorten >=3pt, shorten <=3pt] (vii.west) -- (i.south);

    \end{tikzpicture}
    \caption{Implication cycle in the proof of the Portmanteau Theorem}
    \label{fig:proof-cycle}
\end{figure}
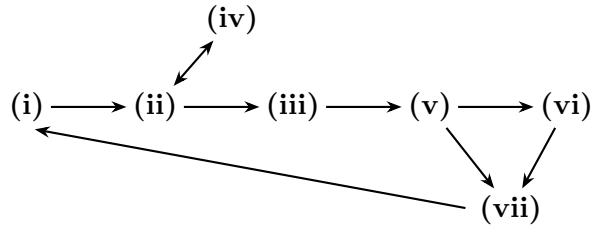

The proof for the Portmanteau theorem proceeds by showing the statements are equivalent via a cyclic chain of implications \citep{Brotherston2022IntroCyclicProofs}. Figure \ref{fig:proof-cycle} provides a visual of the proof. The same order of implications as \citet{vanderVaart1998} is followed.

\begin{Proof}\par
\noindent
\begin{description}[style=nextline,leftmargin=!,labelwidth=4cm]
 \item[(i)$\Rightarrow$(ii)] Since $X_n \overset{d}{\Rightarrow} X$, by Skorokhod's Representation Theorem, there exist random variables $Y_n$ and $Y$ defined on a common probability space such that $Y_n \sim X_n$, $Y \sim X$, and $Y_n \xrightarrow{a.s.} Y$. Because $g$ is continuous, it follows that $g(Y_n) \xrightarrow{a.s.} g(Y)$. Since $g$ is bounded, we apply the Bounded Convergence Theorem to obtain:
$$ \lim_{n \rightarrow \infty} \mathbb{E}[g(X_n)] = \lim_{n \rightarrow \infty} \mathbb{E}[g(Y_n)] = \mathbb{E}[g(Y)] = \mathbb{E}[g(X)] $$
  \item[(ii)$\Rightarrow$(iii)] Noting that the class of bounded continuous Lipschitz\footnote{ibid.} functions $\mathcal{F}_{L}$ are subset to the broader class of bounded continuous functions $\mathcal{F}_{BC}$, i.e. $$\mathcal{F}_L \subset \mathcal{F}_{BC}$$ we can assert $\textnormal{\textbf{(ii)}} \Rightarrow \textnormal{\textbf{(iii)}}$.
  \item[(iii)$\Rightarrow$(v)]  $\forall$ open sets $G$, we can find a sequence of Lipschitz functions $g_m$ such that $$g_m(x) \le g(x) := \ind_{\left\{x \in G\right\}}$$
  and $g_m\nearrow g$ as $m \rightarrow \infty$. For example, $\forall m > 0$ 
  $$
  g_m(x) = \min \left\{1, m\cdot \text{dist}(X, G^c) \right\}, \space \forall x
  $$
  where $\text{dist}(x, G^c) = \underset{y \in G^c}{\inf}d(y,x)$.
  \\
  \\
  So $\forall m >0$
  $$
  \liminf P(X_n \in G) = \liminf \mathbb{E}[g(X_n)] \ge \liminf_{m \rightarrow \infty} \mathbb{E}[g_m(X_n)] \overset{\textbf{(iii)}}{=}\mathbb{E}[g_m(X)] 
  $$
  Since $g_m \nearrow g$, $g_m \le g_{m+1}$. So by the Monotone Convergence Theorem, taking $m \rightarrow \infty$ gives $$P(X_n \in G) \ge P(X \in G) \quad (v)$$
  So we can assert that $\textnormal{\textbf{(iii)}} \Rightarrow \textnormal{\textbf{(v)}}$
  \item[(v)$\Rightarrow$(vi)] Can be shown by taking the compliment of $\liminf P(X_n \in G)$. Noting that the compliment of all open sets $G$ is all closed sets $K$, We have:
  $$
  (\liminf P(X_n \in G))^c = \limsup (P(X_n \in G))^c = \limsup P(X_n \in K) \le P(X\in K)
  $$
  \item[(v)+(vi)$\Rightarrow$(vii)] Denote $\bar{A}$ and $A^\circ$ denote the closure and interior of $A$, then by the fact that $P(X \in \partial A) = 0$
  $$
  P(X \in A) = P(X \in \bar{A}) = P(X \in A^\circ)
  $$
  Then the proof follows by: 
  $$
  \liminf_{n\to \infty} P (X_n \in A) \ge \liminf_{n\to \infty} P(X_n \in A^\circ) \overset{\textnormal{\textbf{(v)}}}{\ge} P(X \in A^\circ) = P(X\in A) 
  $$
  And that:
  $$
  \limsup_{n\to \infty} P(X_n \in A) \le \limsup_{n\to \infty} P(X_n \in \bar{A}) \overset{\textnormal{\textbf{(vi)}}}{\le} P(X \in \bar{A}) = P(X \in A)
  $$
  \item[(vii)$\Rightarrow$(i)] Fixing $x$ as an arbitrary continuity point of $P(X \le x)$. Consider $A = (x, +\infty)$. Then $$P(X \in \partial A) = P(X=x) = 0$$
  by $\textnormal{\textbf{(vii)}}$ we have
  $$
  P(X_n \le x)  \rightarrow P(X \le x) \Longrightarrow X_n \overset{d}{\Rightarrow} X \quad \textnormal{\textbf{(i)}}
  $$
  \item[(ii)$\Leftrightarrow$(iv)] This part is left as an exercise in \citet{vanderVaart1998}. However, we will attempt to solve it here based on what we learned thus far.
  \\
  \\
Let $g$ be a non-negative continuous function. For any constant $M > 0$, we define the truncated function $g_M(x) = \min\{g(x), M\}$. Since $g_M$ is bounded and continuous, assumption \textnormal{\textbf{(ii)}} implies:
$$
\lim_{n\to\infty}\mathbb{E}[g_M(X_n)] = \mathbb{E}[g_M(X)]
$$
Noting that $g(x) \ge g_M(x)$ for all $x$, we have $\mathbb{E}[g(X_n)] \ge \mathbb{E}[g_M(X_n)]$. Therefore:
$$
\liminf_{n\to\infty}\mathbb{E}[g(X_n)] \ge \lim_{n\to\infty}\mathbb{E}[g_M(X_n)] = \mathbb{E}[g_M(X)]
$$
Since $g_M(X)$ is non-decreasing and converges to $g(X)$ as $M \to \infty$, by the Monotone Convergence Theorem we obtain:
$$
\lim_{M \to \infty} \mathbb{E}[g_M(X)] = \mathbb{E}[g(X)]
$$
Combining these results gives the desired inequality:
$$
\liminf_{n\to\infty}\mathbb{E}[g(X_n)] \ge \mathbb{E}[g(X)]
$$
\end{description}
\end{Proof}

\section{Discussion}
The primary objective of this article is to demonstrate that the Portmanteau Theorem can be proven strictly from a probabilistic perspective. While \citet{vanderVaart1998} proves $\mathbf{(i)}\Rightarrow\mathbf{(ii)}$ using a rigorous analytical approach by approximating continuous functions via rectangles and $\varepsilon$-bounds, this proof leverages Skorokhod’s Representation Theorem and the Bounded Convergence Theorem to bypass these complexities.

Comparing this seven-statement framework to the more concise version in \citet{Billingsley1999} and \citet{durrett:2019} highlights a pedagogical trade-off. Billingsley and Durrett prioritizes the topological duality between open and closed sets, providing a streamlined path to the boundary condition ($\mathbf{vii}$). Similarly, \citet{Pollard_2001} maps the cycle of equivalences for weak convergence, emphasizing the connection between functional and topological characterizations. In contrast, by explicitly including intermediate steps such as the Lipschitz class of functions and non-negative expectations, this proof provides a granular view of the logical dependencies. The cyclic chain illustrates how the convergence of distribution functions is fundamentally equivalent to the convergence of expectations across a hierarchy of function classes.

For historical context, it is worth noting that the title ``Portmanteau Theorem" (``portmanteau" translates loosely to ``a large travelling bag or suitcase" in french \citep{mathematics_se}) was introduced by Billingsley in 1968 to describe the utility of packing multiple equivalent definitions of weak convergence into a single theorem. However, the equivalences themselves predate this terminology. The foundational concept of weak convergence on the real line was first developed by Paul L\'evy in the 1920s \citep{levy1925calcul}. The fundamental connection between weak convergence and the preservation of measure in open and closed sets was established by \citet{Alexandroff1943} in the context of general topological spaces. The topological foundations were later integrated into modern probability theory by \citet{Prokhorov1956} and \citet{Skorokhod1956}, whose development of limit theorems for random processes made these equivalences essential to the field. Ultimately, the theorem remains true to its name: a elegantly packed collection of decades of mathematical insight that continues to underpin modern probability theory.
\bibliographystyle{plainnat}
\bibliography{references}

\end{document}